\documentclass[10pt]{article}
\usepackage{amssymb,amsfonts,amsmath,latexsym}
\usepackage{amscd} 
\usepackage{color,float,fancyhdr}
\usepackage[colorlinks=true,linkcolor=blue,citecolor=blue,urlcolor=blue]{hyperref}
\usepackage{graphicx}

\usepackage{soul}

\newtheorem{theo}{\sc Theorem}[section]
\newtheorem{prop}{\sc Proposition}[section]

\newtheorem{cor}{\sc Corollary}[section]
\newtheorem{lem}{\sc Lemma}[section]

\newtheorem{rk}{\sc Remark}[section]

\def\real{\mathbb{R}}

\newcommand{\1}{\mbox{\rm 1\hspace{-0.3em}I}}

\def\k{{\kappa}}

\def\ve{{\varepsilon}}
\def\l{{\lambda}}
\def\d{{\delta}}

\def\un{\underline{n}}
\def\um{\underline{m}}
\def\uy{\underline{y}}
\def\ug{\underline{\gamma}}
\def\uy{\underline{y}}

\def\ue{\underline{e}}

\def\definedas{\stackrel{\Delta}{=}}

\def\grad{\nabla}

\def\HH{{\cal H}}
\newcommand{\LL}{{\cal L}}

\def\NN{{\cal N}}
\def\PP{{\cal P}}
\def\TT{{\cal T}}

\def\qed{$\Box$}

\def\Gr{\text{Graff}}

\def\cov{\text{cov}}
\def\var{\text{var}}

\begin{document}

\title{CLT for Lipschitz-Killing curvatures of excursion sets of Gaussian fields}
  
\author{Marie Kratz\thanks{ESSEC Business School, CREAR, Paris, France; \, E-mail: kratz@essec.edu } \; and  \; Sreekar Vadlamani \thanks{TIFR-Center for Applicable Mathematics (CAM), Bangalore, India; \,  E-mail: sreekar@tifrbng.res.in } }


\maketitle

\begin{abstract}
\noindent Our interest in this paper is to explore limit theorems for various geometric functionals of excursion sets of 
isotropic Gaussian random fields. In the past, limit theorems have been proven for various geometric functionals of
excursion sets/sojourn times ( see \cite{Berman-book, KL-extr, KL-jotp, mesh-shash, Pham, spodarev} for a sample
of works in such settings). The most recent addition being \cite{EL13} where a CLT for Euler-Poincar\'e characteristic 
of the excursions set of a Gaussian random field is proven under appropriate conditions.

In this paper, we shall obtain a central limit theorem for some global geometric functionals, called the Lipschitz-Killing
curvatures of excursion sets of Gaussian random fields in an appropriate setting. \\


\end{abstract}

\section{Introduction and main result}
\label{sec:intro}

There has been a recent surge in interest in understanding the geometry of random sets. In particular,
there have been many works on limit theorems of geometric functionals of random sets coming from
discrete type models arising from various point processes (see \cite{yogi16}, and references therein), or from models of smooth random fields 
\cite{MarCam16}, \cite{EL13}, \cite{probasurvey}, \cite{KL-jotp}, \cite{MarVad16}, \cite{Pham}, \cite{spodarev}.

The object of this paper is to go further, and provide asymptotic distributions for some global 
geometric characteristics of the excursion sets of random fields as the parameter space is allowed to grow to {\it infinity}.

More precisely, let $f$ be a random field defined on $\real^d$, and let $\TT$ be a $d$-dimensional box $[-T,T]^d$ 
We shall be considering the restriction of $f$ to the subset $\TT$, and accordingly 
define the excursion set of $f$ over a threshold $u$, denoted by $A_u(f;\TT)$, as
\begin{equation}\label{Au}
A_u(f;\TT)=  \{ x\in \TT:  ~ f(x)\geq u\} 
\end{equation}
Our interest, in this paper, is to study the distributional aspects of {\it Lipschitz-Killing curvatures} of the sets $A_u(f;\TT)$. 

The Lipschitz-Killing curvatures (LKCs) of a $d$-dimensional {\it Whitney stratified manifold}\footnote{For more details on this, 
and more, we refer the reader to \cite{RFG}.} $M$ are $(d+1)$ integral geometric functionals 
$\{\LL_k(M)\}_{k=0}^d$, with $\LL_0(M)$ the Euler-Poincar\'e characteristic of the set $M$, and 
$\LL_d(M)$ the $d$-dimensional Hausdorff measure of $M$. Though, for $k=1,\ldots, d-1$, the $\LL_k(M)$ do 
not have such clear interpretation, the scaling property\footnote{For any $\lambda>0$, we have $\LL_k(\lambda M)
= \lambda^k \LL_k(M)$, where $\lambda M = \{\lambda x: x\in M\}$.} of the LKCs can be used to interpret the
$k$-th LKC as a $k$-dimensional measure. This property not only underlines the importance of the LKCs, but also
characterises them with their additive and scaling property together with the rigid motion invariance.

One of the most important results in convex geometry is the {\it Hadwiger's characterisation theorem}, 
which, in simple terms, states that the LKCs form a basis for all finitely additive, monotone and rigid motion 
invariant valuations defined on the the collection of {\it basic complexes} (see \cite{RFG}, \cite{KlainRota}),
underlining the importance of LKCs in the study of global geometric characteristics of {\it nice} sets.

The main result of this paper provides a CLT for LKCs of excursion sets of suitable Gaussian random fields.
\begin{theo}
\label{clt-LLk-R}
Let $\TT$ be as defined above, and $f$ be a mean zero, unit variance, isotropic Gaussian random field defined on 
$\real^d$ with $C^3$ trajectories. Then, under some standard regularity assumptions on $f$ as stated in (H1), (H2) and (H3), we have 
\begin{equation}
\frac{ \LL_k\left( A_u(f;\TT)\right) - \mathbb{E}\left( \LL_k\left( A_u(f;\TT)\right)\right)}{|\TT|^{1/2}} \to N(0,\sigma_k^2(u)),\,\,\,\text{as } \TT\to\real^d ,
\end{equation}
for $k=0,\ldots,d$, where $|\TT|$ denotes the $d$-dimensional volume of $\TT$ and,  by $\TT\to\real^d$, we mean $T\to\infty$.
\end{theo}

\begin{rk}
We note here that the specific case corresponding to $k=d$ has already been studied in \cite{Pham}, which is a 
generalization to higher dimensional setting of known results about limit theorems for sojourn times of Gaussian processes.
Another interesting case, namely the Euler-Poincar\'e characteristic (case $k=0$), has been studied in \cite{EL13}, 
whereas a more general result, a CLT for the Euler integral, was obtained in \cite{AdlerNaitzat}. 
\end{rk}

We shall adopt the now standard approach of projecting Gaussian functionals of interest onto the It\^o-Wiener chaos, 
then use the Breuer-Major type of theorem to conclude our main result. This approach has been developed in 
\cite{KL-jotp} to obtain CLT  for general level functionals of $\left( f,\partial f, \partial^2 f\right)$ in dimension $1$, then extended to dimension 2 (see \cite{probasurvey} for a general review on the topic). As applications, they got back CLTs for the number of crossings of $f$, a result first obtained by Slud with an alternative method (see \cite{slud91}, \cite{slud94}), for the number of local maxima (\cite{KL-extr}), for the sojourn time of $f$ in some interval (\cite{probasurvey}), and also for the length of a level curve of a 2-dimensional Gaussian field (\cite{KL-jotp}). 
Note that in these papers, the last step of the method was to approximate $f$ with an $m$-dependent process (\cite{Berman-book}) in order to conclude the CLT. In fact, this step can be removed, and simplified using what is 
now called the {\it Stein-Malliavin method} to conclude Breuer-Major types of theorem, as documented in \cite{noupebook}.
Based on this general approach, CLTs have been proved recently when considering a $d$-dimensional 
Gaussian random field $f$ by Pham \cite{Pham} for the sojourn time of $f$; by Estrade and Le\'on 
\cite{EL13} for the Euler-Poincar\'e characteristic (EPC) of the excursion set of $f$, and by Adler and Naitzat \cite{AdlerNaitzat} for Euler integrals over excursion sets of $f$. Note that EPC shares 
strikingly similar integral representation as the number of level crossings in terms of functional of $f$ but for dimension $d$. Largely, the sketch of the proof for CLTs and the main technical steps remain the same as in \cite{KL-jotp} when dealing with level functionals of a $d$-dimensional Gaussian field $f$, with $d>2$;  the difficulty 
lies in finding a way to avoid explicit computations. Hence the main and crucial contribution in \cite{EL13} has been 
to come up with a neat trick to circumvent this difficulty, proving that the order of the variance of the EPC for $f$ 
restricted to any subspace of $\TT$ is less than $|\TT|$, hence is negligible in the limit  as $\TT$ grows to $\real^d$.

We are going to build on those works, in order to obtain a CLT for all LKCs of excursion sets. The 
difficulty here is to develop similar techniques when working on $A_u(f;\TT\cap V^*)$ for any $k$-dimensional 
affine subspace $V^*$ of $\real^d$. 

The structure of this paper is as follows.
Throughout this paper, we work with isotropic Gaussian random fields. 
We  begin in Section \ref{sec:prelim} with setting the notation, and the necessary background  
for the analysis to follow in later sections. In Section \ref{subsec:chaos}, we recall the basics of
the expansions of Gaussian functionals using the multiple Wiener-It\^o integrals. Next, in Section 
\ref{subsec:LKC-Crofton}, we  define the Lipschitz-Killing curvatures, and also state the {\it Crofton formula} 
that provides a relationship between various LKCs; it is going to be a crucial element in the proof of our 
main result. Section \ref{subsec:EPC-ex-set} is devoted to discuss the integral representation of Euler-Poincar\'e 
characteristic of excursion sets of any random field via the {\it expectation metatheorem}. Finally, precise setup 
for the problems, and the assumptions, in particular on the covariance structure of the random field $f$, is listed 
in Section~\ref{subsec:assume}. In Section~\ref{sec:main-proof}, we develop the proof of the main result, 
Theorem~\ref{clt-LLk-R}, using the standard sketch given in three main steps. First  we prove that the functional of interest is square-integrable (Section~\ref{subsec:sq-int}) and obtain its Hermite expansion in Section~\ref{subsec:Hermite}. 
Then we prove that the limiting variance is bounded away from zero and infinity in Section~\ref{subsec:var-bounds}. 
Finally, in Section~\ref {subsec:BM}, we give an extension of Breuer-Major theorem to affine Grassmannian case  to conclude the Gaussianity of the limiting distribution. 
Section \ref{sec:discussion} concludes with a discussion and a multivariate CLT for EPCs.

\section{Preliminaries}
\label{sec:prelim}

\subsection{It\^o-Wiener chaos expansions}
\label{subsec:chaos}

Hermite polynomial expansions or Multiple Wiener-It\^o Integrals (MWI) are a powerful tool to approximate
and study nonlinear functionals of stationary Gaussian fields (see \cite{EL13}, \cite{probasurvey}, or \cite{KL-jotp} for details). 

Formally, let $Z$ be a $m$-dimensional standard Gaussian random vector, and $L^2(Z)$ be the set of 
all real square integrable functionals of $Z$. In short, giving a Hermite expansion is a way to approximate elements from $L^2(Z)$ by a series of Hermite polynomials. More precisely, for 
$\un\in \left( \mathbb{N} \cup \{0\}\right)^m$, define $\HH_{\un}(Z)= \prod_{i=1}^mH_{n_i}\left(Z_i\right)$, 
and set $\HH_q$ as the linear span of $\{\HH_{\un}(Z) : |\un| = \sum_{i=1}^m n_i =q\}$. Then 
\begin{equation}\label{eqn:Wiener-chaos}
L^2(Z)=\bigoplus_{q\ge 0} \HH_{q} 
\end{equation}
For more details regarding Hermite expansions, and their applications to study functionals
of Gaussian random fields, we refer the reader e.g. to \cite{az-ws}, \cite{Berman-book}, \cite{noupebook}. 

The above can also be written in a more abstract setting of multiple Wiener integrals, for which we begin
with an orthonormal system $\{\varrho_i\}_{i\ge 1}$ of $L^2(\real^m)$. Writing $W$ for complex Brownian measure
on $\real^m$, let us define $\xi_i = \int_{\real^m} \varrho_i(\lambda)\,W(d\lambda)$. 
Clearly, $\{\xi_i\}$ form a sequence of i.i.d. standard normal random variables. Now for a fixed
$\un\in \left(\mathbb{N}\cup \{0\}\right)^q$, and $p_1,\ldots,p_q \in \mathbb{N}$ define 
$H_{\un}(\xi_{p_1},\ldots,\xi_{p_q}) = \prod_{i=1}^q H_{n_i}(\xi_{p_i})$. Then (see \cite{Nualart-book})
\begin{eqnarray}
\label{eqn:multiple-Wiener-def}
H_{\un}(\xi_{p_1},\ldots,\xi_{p_q}) &=& \int_{\real^{mq}} 
\left( \varrho^{\otimes n_1}_{p_1}\otimes\cdots\otimes\varrho^{\otimes n_q}_{p_q}\right)
(\lambda_1,\ldots,\lambda_q) \, W(d\lambda_1)\cdots W(d\lambda_q)\nonumber\\
&\definedas & I_q\left( \varrho^{\otimes n_1}_{p_1}\otimes\cdots\otimes\varrho^{\otimes n_q}_{p_q}\right)
\end{eqnarray}
where $I_q$ denotes multiple Wiener integral. A decomposition, similar to \eqref{eqn:Wiener-chaos},
holds true for all square integrable functionals of $W$, and is called the {\it It\^o-Wiener chaos}. We refer the reader 
to \cite{Nualart-book} for complete details.

\subsection{Lipschitz Killing curvatures and the Crofton formula}
\label{subsec:LKC-Crofton}

There are a number of ways to define Lipschitz-Killing curvatures, but
perhaps the easiest is via the so-called Weyl's tube formula (see \cite{Hotelling}, \cite{Wey39} for the first hand account of this formula). In order to
state the tube formula, let $M$ be an $m$-dimensional manifold with {\it positive reach} (see \cite{RFG}) 
embedded in $\real^{n}$ which is endowed with the canonical Riemannian structure on $\mathbb{R}^{n}$. 
Then, writing $\|\cdot\|$ as the standard Euclidean norm on $\real^n$, the tube of radius $\rho$ 
around $M$ is defined as 
\begin{equation}
\text{Tube}(M,\rho )\ =\ \left\{ x\in \mathbb{R}^{n}:\;\inf_{y\in M}\Vert x-y\Vert \leq \rho \right\}.
\end{equation}%

Then according to Weyl's tube formula {(see \cite{RFG})}, the Lebesgue
volume of so constructed tube, for small enough $\rho$, is given by 
\begin{equation}
\lambda _{n}(\text{Tube}(M,\rho ))\ =\ \sum_{j=0}^{m}\rho ^{n-j}\omega _{n-j}%
\LL_{j}(M)\text{ },  \label{tube:formula}
\end{equation}%
where $\omega _{n-j}$ is the volume of the $(n-j)$-dimensional unit ball in $\real^{n-j}$, and $%
\LL_{j}(M)$ is the $j$-th LKC of $M$.
Although, it may appear from the definition above that the $\LL_{j}$
depend on the embedding of $M$ in $\mathbb{R}^{n}$, in fact, the $\LL_{j}(M)$ are
intrinsic, and so are independent of the ambient space.

Apart from their appearance in the tube formula (\ref{tube:formula}), there
are, at least, two more ways in which to define the LKCs (see \cite{RFG}).

Borrowing the notations from \cite{RFG}, let {\bf $\Gr(d, k)$} be the affine Grassmannian of all $k$-dimensional 
affine subspaces of $\real^d$, and {\bf $\text{Gr}(d,k)$} be the set of all $k$-dimensional linear subspaces of $\real^d$.

Let $M$ be a compact subset of $\real^d$ and $V^*\in\Gr(d,k)$. 
Then writing 
$$
M_{V^*} \;\text{for} \; (M\cap V^*),
$$ 
and setting $\displaystyle \lambda^d_{k}$  to be the appropriate, normalized 
measure on $\Gr(d,k)$ (cf.\cite{RFG}), and also
$$ 
\left[\begin{array}{c} m\\n \end{array}\right] = 
\frac{\omega_m}{\omega_n\,\omega_{m-n}}\left(\begin{array}{c} m \\n \end{array}\right),
$$ 
we have the Crofton formula:
\begin{equation}\label{euclid-crofton}
\int_{\Gr(d, k)} \LL_j(M_{V^*}) \, d\lambda_{k}^d (V^*)= \left[\begin{array}{c} d-k+j\\j \end{array}\right] \LL_{d-k+j}(M)
\end{equation}
whenever $M$ is {\it tame} and a {\it Whitney stratified space} (see \cite{RFG}).

Setting $j=0$ in the above equation \eqref{euclid-crofton} gives back the Hadwiger formula 
\begin{equation}\label{hadwiger}
\int_{\Gr(d, k)} \LL_0(M_{V^*}) d\lambda_{k}^d (V^*) = \LL_{d-k}(M),
\end{equation}
which we shall use to generate all the LKCs given the Euler-Poincar\'e characteristic of 
all the slices $M_{V^*}$.

Another interesting case is when we set $j=k$ in \eqref{euclid-crofton}; we obtain
\begin{equation}\label{crofton-volume}
\int_{\Gr(d, k)} \left|M_{V^*}\right| d\lambda_{k}^d (V^*)= \left[\begin{array}{c} d\\k \end{array}\right] \LL_{d}(M)
\end{equation}
where $\left|M_{V^*}\right|$ is the $k$-dimensional Hausdorff measure of the set $M_{V^*}$.

\subsection{Euler-Poincar\'e characteristic and other LKCs of excursion sets}
\label{subsec:EPC-ex-set}

Let $\TT$ be a compact, tame and Whitney stratified subset of $\real^d$. For any fixed 
$V^*\in\text{Graff}(d,k)$, set $\partial_l \TT_{V^*}$ as the  $l$-dimensional boundary of $\TT_{V^*}$. 
Assume $f$ be a smooth Gaussian random field, then using the standard Morse theory 
(see \cite[Chapter 9]{RFG}), we can write
\begin{equation} \label{eqn:EPC}
\LL_0\left( A_u(f;\TT_{V^*})\right)= \sum_{l=0}^{k} \sum_{J\in \partial_l \TT_{V^*}} \phi_l(J)
\end{equation}
whenever $\TT$ is {\it tame and a Whitney stratified space} (see \cite{RFG}), where,
$$\phi_l(J) 
= \sum_{j=0}^l (-1)^j \# \left\{ t\in J: f(t) \ge u, \ \grad_J f(t) = 0,\ \text{index}(\grad_J^2 f(t))=l-j\right\},$$
with  $\grad_J f$ and $\grad^2_J f$ representing restrictions of the usual gradient $\grad f$ and Hessian $\grad^2 f$ 
onto $J\in\partial_l \TT_{V^*}$. 

Applying Theorem 11.2.3 of \cite{RFG}, the above equation can formally be rewritten as
\begin{equation}\label{eqn:phi-integral}
\phi_l(J) 
=  (-1)^l \int_{\TT_{V^*}} \d(\grad_J f(t)) \,\1_{\{f(t)\ge u\}} \, \det\left( \grad_J^2 f(t)\right) \, dt
\end{equation} 
almost surely and in $L^2$, where $\delta$ is the Dirac delta at $0$ defined on $\real^d$, 
interpreted as usual by approximating $\delta$, as $\ve\to 0$, by the Gaussian density of a 
$d$-vector with independent components mean  $0$ and variance $\ve$, or by the function 
$\displaystyle (2\ve)^{-d}\1_{[-\ve;\ve]^d}$ (see e.g. \cite{RFG} for a.s., and \cite{KL-SPA} 
or \cite{EL13} for $L^2$ convergence). Hence the way of obtaining a Hermite expansion 
of $\LL_{d-k}$ will go through a limiting process. However, this process of approximation 
being clearly spelled out in many of previous works going as far back as \cite{Berman-book}, 
we shall omit this step in the rest of the paper, and skip to the limit.

We shall now combine equations \eqref{hadwiger} and \eqref{eqn:EPC} to express all other LKCs in terms of
the Euler-Poincar\'e characteristic of $A_u(f;\TT_{V^*})$. Formally,
\begin{equation}\label{eqn:Lk-Au}
\LL_{d-k}\left( A_u(f;\TT)\right) = \int_{\text{Graff}(d,k)}\LL_0\left( A_u(f;\TT_{V^*})\right) \,d\lambda^d_{k}(V^*).
\end{equation}

\begin{rk} - {\bf Parametrization of $\text{Graff}(d,k)$}
\label{rk:parametrization-graff}~\\
Note that $\text{Graff}(d,k)$ can be parametrized as
$\text{Gr}(d,k)\times \real^{d-k}$. Furthermore, we shall identify $\text{Gr}(d,k)$ with the set of all $k\times d$ 
matrices whose rows are orthonormal vectors in $\real^d$, modulo left multiplication by a $k\times k$ orthogonal matrix (see\cite{Rubin}).
\end{rk}

Writing $V$ as the matrix whose rows are $k$-orthonormal vectors spanning the linear space 
obtained by the parallel translate of $V^*$, 
\begin{equation} \label{eqn:critical-pts}
\phi_k(\partial_{k} \TT_{V^*}) 
= (-1)^{k} \int_{\TT_{V^*}} \d(V\grad f(t)) \, \1_{\{f(t)\ge u\}} \, \det\left( V\grad^2 f(t) V^T\right) \, dt.
\end{equation} 

\noindent {\sc Note:} For any {\bf $V^*\in\Gr(d,k)$}, we shall denote $V$ for the matrix
whose rows are $k$-orthonormal vectors spanning the linear space 
obtained by the parallel translate of $V^*$, and we shall use the same $V$ to denote the 
element in $\text{Gr}(d,k)$ that corresponds to the $k$ dimensional linear space
spanned by the rows of the matrix $V$.

\subsection{Setup for the problem and assumptions}
\label{subsec:assume}

In this paper, we consider $f$ a mean zero, isotropic, real valued Gaussian random field defined on $\real^d$ with $C^3$ trajectories.
The assumption of isotropy means that the covariance of the Gaussian random field satisfies
$$\mathbb{E}[f(x) f(y)]= r(|x-y|), \,\,\,\,\,\,\,\,\forall x,y\in\real^d$$
for some function $r:\real_{+} \to \real$. Without loss of generality, we shall assume $r(0)=1$.

We denote the partial derivatives of order $n$ of any function $g$ defined on $\real^d$ as
$$
g^{(i_1 \cdots i_n)}(t)=\frac{\partial ^n}{\partial t_{i_1}\cdots\partial t_{i_n}} g(t).
$$

We introduce the gradient $\grad f(x)$ and Hessian $\grad^2 f(x)$ of $f$, and recall that due to isotropy 
$\grad f(x)$ and $(f(x), \grad^2 f(x) )$ are independent for every fixed $x$ (in fact, stationarity suffices to conclude the same).  
Thus, the covariance function of $(\grad f(x), f(x), \grad^2 f(x))$ can be expressed as a block diagonal matrix for each fixed $x$. 
We denote the covariance matrix of $(\grad f(x), f(x), \grad^2 f(x))$ as
\begin{equation}\label{eqn:Lambda}
\Sigma = \left(
\begin{array}{cc}
\Sigma_1 & 0 \\ 0 & \Sigma_2
\end{array}
\right)
\end{equation}
where $\Sigma_1$ and $\Sigma_2$ are the covariance matrices of $\grad f$, and $(f, \grad^2 f)$, respectively.
Notice that since $(\grad f(x), f(x), \grad^2 f(x))$ is a $\left(d + 1 + d(d+1)/2\right)$ dimensional vector, the corresponding
covariance matrix is a square matrix of order $D\times D$, where 
\begin{equation}\label{def:D}
D=d + 1 + d(d+1)/2.
\end{equation}
Simple linear algebraic considerations imply that there exists a $D\times D$ matrix $\Lambda$ such that $\Sigma =\Lambda\Lambda^T $. Then we define a new field  $Z=(Z(x), x\in\real^d)$ by
\begin{equation}\label{Z} 
Z(x)= \Lambda^{-1} \left(\grad f (x), \grad^2 f (x), f (x) \right).
\end{equation}
Let us denote its covariance function by
\begin{equation}\label{eqn:covZ}
\gamma=\left(\gamma_{ij}(.)\right)_{1\le i,j\le D} \quad \text{with}\; \,\gamma_{ij}(h)= \cov\left(Z_i(x),Z_j(x+h)\right).
\end{equation}
Note here that we have implicitly used the fact that various derivatives of a stationary Gaussian random field are themselves stationary Gaussian random fields (see \cite[Chapter 5]{RFG}).

We can also write $Z$ as $\displaystyle Z(x)=  (Z^{(1)}(x),Z^{(2)}(x))$ so that
\begin{equation}\label{eqn:Z1-Z2}
Z(x) = (Z^{(1)}(x),Z^{(2)}(x))\sim \NN(\underline{0},I_D), \quad I_D \;\text{being  the identity in} \, \real^D.
\end{equation} 

We need more assumptions on $f$ to ensure that various LKCs of the excursion set $A_u(f; \TT)$ of $f$ over a threshold $u$,  are indeed square integrable, and that they satisfy a CLT as $\TT\to\real^d$. 

The required assumptions are rather standard when looking for CLT of non linear functionals of stationary Gaussian random fields, such as number of crossings \cite{probasurvey,KL-jotp},  curve length \cite{KL-jotp}, EPC \cite{AdlerNaitzat,EL13}, sojourn time \cite{Pham}, etc.

\begin{itemize}
\item[(H1)] {\it Geman type condition:
We shall assume that the covariance function $r\in C^4 (\TT)$, and that
the function $\displaystyle \frac1{\|t\|^{2}}\left( \grad^2 r(t) -r^{(ii)}(0) I_d\right)$ defined on $\real^d$, is bounded near 
$t=0$, where we recall that $r^{(ii)}(0)=-\var(\grad_i f(x))$. }\\
\vskip0.1cm
(H1) is simply the higher dimensional analog of Geman's condition (\cite{geman}), which is needed to prove that 
the functional of interest is in $L^2$. When $d=1$, it is known to be a necessary and sufficient condition to obtain the $L^2$ convergence of the number of crossings of any threshold (see \cite{KL-AP}).
Observe that this condition (H1) is satisfied whenever the underlying random field is `smooth' enough, with   $C^3$ sample paths\footnote{We refer
to \cite{az-ws} for further discussion on connections between analytical assumptions and regularity conditions for Gaussian random fields.}. 
Hence, for simplicity, we will assume now on that $f$ {\bf is} $\mathbf C^3$. 
\vskip0.2cm

\item[(H2)] {\it Arcones type condition:}
{\it For the covariance function $\gamma(\cdot)$ of the joint field $\left(Z(x)\right)_{x\in \real^d}$,
we assume that there exists an integrable $\psi$ on $\real^d$ satisfying 
$$\psi(t) \underset{||t||\to +\infty}{\longrightarrow} 0,$$ such that
$$
\max_{1\le i,j \le D} |\left(\gamma(x)\right)_{ij}| \le K \psi(x), \quad \text{ for some }K > 0.
$$
}
\vskip0.1cm
This condition is crucial in ensuring the finiteness of the limiting variance of the considered functionals.
As already noted in \cite{KL-jotp} and \cite{EL13}, it implies in particular the existence of the spectral density, and that 
$r \in L^q(\real^d)$, $q\ge 1$.
\vskip0.3cm

\item[(H3)] {\it The spectral density, denoted by $h$, of the covariance function corresponding to the field $f$ satisfies $h(0) >0$. }\\
\vskip0.1cm
This condition is needed  to ensure that the asymptotic variance obtained in the CLT is non zero.

\end{itemize}

\section{Proof of Theorem \ref{clt-LLk-R}}
\label{sec:main-proof}

Recall that the key argument in \cite{EL13} to prove the CLT for the Euler-Poincar\'e
characteristics of the excursion set is to consider only the highest dimensional term in \eqref{eqn:EPC}, dropping all lower dimensional terms, and proving later that the contribution from the lower dimensional terms is negligible under the volume scaling. We will also use the same argument. 

Let us begin with 
\begin{eqnarray} \label{eqn:L0-sum-dim}
&& \LL_{d-k}\left( A_u(f;\TT)\right) \nonumber\\
& = & \int_{\text{Graff}(d,k)}\sum_{l=0}^{k} (-1)^l \sum_{J\in\partial_l \TT_V} 
\int_{J}  \1_{(f(x) \ge u)}\,  \d (\grad_J f(x)) \det\left(\grad_J^2 f(x)\right) dx \,\,d\lambda^d_{k}(V) \nonumber\\
&=& \sum_{l=0}^{k} \LL_{d-k,l}\left( A_u(f;\TT)\right)
\end{eqnarray}
where 
\begin{eqnarray*}
&& \LL_{d-k,l}\left( A_u(f;\TT)\right)\\
&=& (-1)^l \int_{\text{Graff}(d,k)} \sum_{J\in\partial_l \TT_V} 
\int_{J}  \1_{(f(x) \ge u)}\,  \d (\grad_J f(x)) \det\left(\grad_J^2 f(x)\right) dx \,\,d\lambda^d_{k}(V).
\end{eqnarray*}

We shall consider only $\LL_{d-k,k}\left( A_u(f;\TT)\right)$, and prove that, after appropriate normalization, it exhibits a central limit theorem. 
Thereafter, the same arguments can be pieced together to conclude that, under the same scaling, $\LL_{d-k,l}\left( A_u(f;\TT)\right)$ converge to $0$ whenever $l<k$.

As spelt out in the introduction, our proof has three major steps:
\begin{enumerate}
\item to show that the functional of interest is square integrable, and thereby obtain its Hermite type expansion;
\item to prove that the limiting variance is bounded away from zero and infinity;
\item to use the Stein-Malliavin method for Breuer-Major type functionals to conclude to the 
Gaussianity of the limiting distribution.
\end{enumerate}

We shall provide details of the aforementioned steps in the following subsections.

\subsection{Square integrability}
\label{subsec:sq-int}

Let us recall \eqref{eqn:critical-pts}, and write
\begin{equation}\label{eqn:def-LL(d,d-k)}
\LL_{d-k,k}\left(A_u(f;\TT)\right) = \int_{\Gr(d,k)} \phi_k(\partial_k \TT_{V^*}) \, d\lambda^d_k(V^*).
\end{equation}
Clearly, $\phi_k(\partial_k \TT_{V^*}) = 0$ whenever $\TT_{V^*}$ is empty. Moreover, compactness of $\TT$
implies that the domain of integration in the above integral is a compact bounded subset of $\Gr(d,k)$ given
by ${\cal W} = \{V^*\in\Gr(d,k): \TT_{V^*} \text{ is non empty} \}$. Also, under the standard normalization of
$\lambda^d_k$ (cf. \cite{RFG}), we have $$\lambda^d_k({\cal W}) = |\TT| \left[\begin{array}{c} d\\ k\end{array}\right]$$
Applying Jensen's inequality, we have
\begin{equation} \label{eqn:second-moment}
\mathbb{E}\left[\left(\LL_{d-k,k}\left(A_u(f;\TT)\right)\right)^2\right] 
\le |\TT| \left[\begin{array}{c} d\\ k\end{array}\right] \int_{\Gr(d,k)} \mathbb{E}\left[\left(\phi_k(\partial_k \TT_{V^*})\right)^2\right] \, d\lambda^d_k(V^*).
\end{equation}
We shall now focus on obtaining an appropriate upper bound for the integrand in the above expression, 
which in turn shall imply the square integrability of $\LL_{d-k,k}\left(A_u(f;\TT)\right)$.

Using standard fare, notice that $\phi_k(\partial_k \TT_{V^*})$ can be bounded above by the cardinality of the set $\{t: V\grad f=0\}$ 
that we denote by $N_u(\TT_{V^*})$. Then, as usual, we compute the
second factorial moment of $N_u(\TT)$ and prove that it is finite to conclude the square integrability.
We have, setting $\Xi_V(t_1,t_2)$ as the set $\{V\grad f(t_1) = V\grad f(t_2) = 0\}$, and using \cite[Corollary 11.5.2]{RFG}, or \cite[Theorem 6.2]{az-ws}
\begin{eqnarray}\label{eqn:L2-factorial}
&& \mathbb{E}\left[ N_u(\TT_{V^*})\left(N_u(\TT_{V^*})-1\right)\right]  \,= \nonumber\\
& = & \int_{\TT_{V^*}^2} \!\!\!\!\!
\mathbb{E}\left(\left. |\det(V\grad^2 f(t_1)V^T)\det(V\grad^2 f(t_2)V^T)|\,\right|\Xi_V(t_1,t_2)\right) \nonumber\\ 
&& \hskip 1.5cm\times \,\,\,\,\,\, p_{V,t_1,t_2}(0,0)\,dt_1 dt_2\qquad
\end{eqnarray}
where $\TT^2_{V^*} = 
\{(t_1,t_2)\in \TT_{V^*}\times \TT_{V^*} : t_1\ne t_2\}$, and $p_{V,t_1,t_2}(0,0)$
is the joint density of $\left( V\grad f(t_1),V\grad f(t_2)\right)$. Using stationarity,  we can
reduce the above integral to
\begin{eqnarray}\label{eqn:L2-factorial-1}
&  & \mathbb{E}\left[ N_u(\TT_{V^*})\left(N_u(\TT_{V^*})-1\right)\right] \nonumber\\
&& \!\!\!\!\!\!\!\!\!\!\!\!\!\! = \!\!\!\!\!\!\! \underset{\scriptsize\begin{array}{c}t_1\in \TT_{V^*}\\ s\in (\TT_{V^*}-t_1)/\{0\}\end{array}}{\int} \!\!\!\!\!\!\!\!\!\!\!\!\!\!\!\!\!
\mathbb{E}\left(\left. |\det(V\grad^2 f(0)V^T)\det(V\grad^2 f(s)V^T)| \,\right| \Xi_V(0,s\right) p_{V,0,s}(0,0)\,ds \, dt_1. \nonumber
\end{eqnarray}
Next, using the Cauchy-Schwarz inequality and stationarity gives
\begin{eqnarray*}
&& \mathbb{E}\left[\left. |\det(V\grad^2 f(0)V^T)\det(V\grad^2 f(s)V^T)|\,\right| \Xi_V(0,s)\right]\\
&\le & \mathbb{E}\left[\left. |\det(V\grad^2 f(0)V^T)|^2\right| \Xi_V(0,s)\right].
\end{eqnarray*}
Invoking similar delicate analysis as in the appendix of \cite{EL13}, we can conclude that there exists a constant $C_1$ 
(independent of $T_{V^*}$) such that
\begin{equation}\label{eqn:determinant-bound}
\mathbb{E}\left[\left. |\det(V\grad^2 f(0)V^T)|^2\right| \Xi_V(0,s)\right] \le C_1 \|s\|^2.
\end{equation}
Next, notice that 
\begin{equation}\label{eqn:density-bound}
p_{V,0,s}(0,0) \le  C_2 \|s\|^{-k}
\end{equation}
where $C_2$ is a constant independent of $V$.
Then, combining equations \eqref{eqn:determinant-bound} and \eqref{eqn:density-bound} in \eqref{eqn:L2-factorial} provides
\begin{equation}\label{eqn:L2-factorial-2}
\mathbb{E}\left[ N_u(\TT_{V^*})\left(N_u(\TT_{V^*})-1\right)\right] \le C_3
\!\!\!\!\!\!\!\!\!\!\!\!\!\! \underset{\scriptsize\begin{array}{c}t_1\in \TT_{V^*}\\ 
s\in (\TT_{V^*}-t_1)/\{0\}\end{array}}{\int} \!\!\!\!\!\!\!\! \|s\|^{2-k}\,ds\, dt_1
\end{equation}
for some finite positive constant $C_3$.
Notice first that $(\TT_{V^*}-t_1) \subset 2\TT_{V}$, where $V$
is the translate of $V^*$ containing origin, implying
$$\underset{\scriptsize\begin{array}{c}t_1\in \TT_{V^*}\\ 
s\in (\TT_{V^*}-t_1)/\{0\}\end{array}}{\int} \!\!\!\!\!\!\!\! \|s\|^{2-k}\,ds\, dt_1 \le
\underset{\scriptsize\begin{array}{c}t_1\in \TT_{V^*}\\ 
s\in 2\TT_{V}/\{0\}\end{array}}{\int} \!\!\!\!\!\!\!\! \|s\|^{2-k}\,ds\, dt_1 
= \left| \TT_{V^*}\right| \int_{s\in 2\TT_{V}/\{0\}} \|s\|^{2-k}\,ds.$$
Further, the above integral can be bounded from above by replacing the domain of integration by
the set $[B_d(0,c_d T)\cap \TT_{V}]/\{0\}$, where $B_d(0,r)$ represents a $d$-dimensional ball
of radius $r$, and the constant $c_d$ can be chosen appropriately so as to encompass the set $\TT_V$.
Then, observe that the latter simplified integral does not depend on the choice of $V$. Therefore choosing
$V$ as the span of any $k$ coordinate axes, we obtain
$$\int_{s\in 2\TT_{V}/\{0\}} \|s\|^{2-k}\,ds \le \int_{B_d(0,c_d T)/\{0\}} \|s\|^{2-k}\,ds = \kappa \int_0^{c_d T} r\,dr$$
where the equality is a result of a simple polar coordinate transformation and $\kappa$ is appropriate universal constant.
Therefore, 
\begin{equation}\label{eqn:L2-factorial-3}
\mathbb{E}\left[ N_u(\TT_{V^*})\left(N_u(\TT_{V^*})-1\right)\right] \le C_4 |\TT|^{2/d} \left| \TT_{V^*}\right|
\end{equation}
for some positive constant $C_4$ independent of the choice of $V^*$.

Finally, using the standard Gaussian kinematic fundamental formula (see \cite{RFG}) for the mean 
of $N_u(\TT_{V^*})$ and the above
computations together with equation \eqref{eqn:second-moment} and Crofton formula, we can conclude that
$$\mathbb{E}\left( \LL_{d-k,k}\left(A_u(f;\TT)\right)\right)^2 < C \left|\TT\right|^{2(1+1/d)}$$
for some large, but finite and positive constant $C$.
Note that this upper bound is not optimal, but still suffices to achieve the goal
of square integrability to obtain a Hermite type expansion of the functionals of interest. Using the Hermite type expansion, we shall obtain much tighter bounds later in Section \ref{subsec:var-bounds}.

\subsection{Hermite expansion}
\label{subsec:Hermite}

Set  $F(x) = \left(f_1(x), f_2(x)\right)$ with $f_1(x) = \grad f(x)$ and
$f_2(x) = \left( \grad^2 f(x), f(x)\right)$.

For $x\in \TT$, recall (see Section \ref{subsec:assume}) that we can factorize $\Sigma$ as 
$\Sigma = \Lambda \Lambda^T$, such that
$\Lambda$ has a block diagonal form
\begin{equation}\label{eqn:lambda}
\Lambda = \left(
\begin{array}{cc}
\Lambda_1 & 0 \\ 0 & \Lambda_2
\end{array}
\right)
\end{equation}
where $\Lambda_1$ is the formal {\it square root} of $\Sigma_1$ and $\Lambda_2$ is a lower triangular matrix 
such that $\Lambda_2\Lambda^T_2=\Sigma_2$, respectively.

Using equation \eqref{eqn:critical-pts} and standard methods as in \cite{EL13}, we now obtain 
a Hermite expansion for $\phi_k(\partial_k \TT_{V^*})$. Define
\begin{align}
G_{1}^V(f_1(x)) =&  \d(V\grad f(x))   \label{eqn:G1-def}\\
G^V_{2,u}(f_2(x)) =& \1_{(f(x) \ge u)}\, \det\left(V\grad^2 f(x)V^T\right) .  \label{eqn:G2-def}
\end{align}
Clearly, for each fixed space point $x$, the functions $G_1$ and $G_2$ are independent. We shall obtain Hermite expansions for these two functions separately. 

Formally, for $\un\in\mathbb{N}^D$, $D$ being defined in \eqref{def:D}, set 
$\un = (\un_1, \un_2) \in \mathbb{N}^d\times \mathbb{N}^{D-d}$, then the square integrability
implies that we have the following
\begin{eqnarray}
\label{chaosGeps}
&& G^V_{1}(f_1(x)) \times G^V_{2,u}(f_2(x)) \nonumber\\
& & = \sum_{q=0}^{\infty} \sum_{\un^{(D)}:\sum_{i=1}^D n_i =q} c(\un,u,V,\Lambda)  
H_{\un_1}(Z^{(1)}(x))\, H_{\un_2}(Z^{(2)}(x))
\end{eqnarray}
where the Hermite coefficients are given by
\begin{eqnarray*}
c(\un,u,V,\Lambda) & := & \frac1{\un !} \int_{\real^{D}} G^V_{u}(\Lambda\uy) 
\prod_{i=1}^{D} H_{n_i}(y_i)\,\varphi_{D}(\uy)\,d\uy\\
& = & c_1(\un_1, V, \Lambda_1)\times c_2(\un_2,u,V, \Lambda_2)
\end{eqnarray*}
writing $\varphi_D$ for the standard normal density  in $D$-dimensions.
\begin{rk}\label{rem:invariance}
It follows from the discussion of \cite[Section 5.7]{RFG}, that the distribution of $\left( V\grad f(x), V\grad^2 f(x)V^T\right)$
does not depend on the space point $x$ (due to stationarity) and on $V$ (due to isotropy). Therefore, the coefficients
$c_1(\un_1, V, \Lambda_1)$ and $c_2(\un_2, u, V, \Lambda_2)$ do not depend on $V$, which will help simplifying 
the proofs.
\end{rk}

\noindent $\underline{\text{{\scshape Computing $c_1$}}}$

Observe that
$$
c_1(\un_1, V,\Lambda_1) = \frac1{\un_1!}\int_{\real^d} G^V_{1}(\Lambda_1\uy_1) H_{\un_1}(\uy_1)\,
\varphi_d(\uy_1)\,d\uy_1.
$$
First, note that the integral is to be interpreted as a limit of integral of an appropriate approximation of 
$\delta$. Secondly, by Remark \ref{rem:invariance}, we can choose a $V$ which suits our purpose.

In particular, one may define
$$G^V_{1,\ve}(\Lambda_1\uy_1) = \frac{1}{(2\pi\ve^2)^{k/2}} 
\exp\left( -\frac1{2 \ve^2} \uy_1^T\Lambda^T_1V^TV\Lambda_1\uy_1 \right),$$
and thus define $c_1(\un_1, V,\Lambda_1)$ as an $L^2$-limit  of
$$
c_1(\un_1,\ve, V,\Lambda_1) \definedas \frac1{\un_1! (2\pi\ve^2)^{k/2}}\int_{\real^d}  
\exp\left( -\frac1{2 \ve^2} \uy_1^T\Lambda^T_1V^TV\Lambda_1\uy_1 \right)\, H_{\un_1}(\uy_1)\,
\varphi_d(\uy_1)\,d\uy_1 .
$$
Noticing that the variance of $\grad_i f(x)$ does not depend on the index $i$ due to {\it isotropy}, 
we conclude that $\Lambda_1 = \sqrt{\lambda}\,I_d$ where $\lambda$ is the variance of $\grad_i f(x)$.
Therefore,
$$
c_1(\un_1,\ve, V,\Lambda_1) = \frac{1}{\un_1! (2\pi\ve^2)^{k/2}}\int_{\real^d} 
\exp\left( -\frac1{2 \ve^2} \lambda \,\uy_1^TV^TV\uy_1 \right)\, H_{\un_1}(\uy_1)\, \varphi_d(\uy_1)\,d\uy_1 .
$$
Equivalently,
$$
c_1(\un_1,\ve, V,\Lambda_1) = \frac{\lambda^{-k/2}}{\un_1!}\int_{\real^d} 
\frac{\lambda^{k/2}}{\ve^k}\varphi_k\left(\frac1{\ve}\sqrt{\lambda} V\uy_1\right)\, H_{\un_1}(\uy_1)\, \varphi_d(\uy_1)\,d\uy_1,
$$
where $\frac{\lambda^{k/2}}{\ve^k}\varphi_k(\sqrt{\lambda}\ve V\uy_1)$ converges to the desired Dirac delta. 
As pointed earlier, the above computation is invariant of the choice of $V$, so we shall choose $V$ to be the space spanned by
$(e_1,\ldots,e_k)$ where $\{e_i\}_{i=1}^d$ is the canonical basis of $\real^d$. Thereafter, taking limit as $\ve\to 0$,
we obtain
\begin{equation}\label{eqn:c1-precise}
c_1(\un_1, V,\Lambda_1) = (2\pi\lambda)^{-k/2}\prod_{i=1}^{k}\frac{H_{n_{1,i}}(0)}{n_{1,i}!}.
\end{equation}
However, in order to obtain estimates for the limiting variance, we shall need bounds on $c_1(\un_1, V,\Lambda_1)$.
Using the usual technique as sketched in \cite{KL-SPA}, we obtain
\begin{eqnarray*}
\left|c_1(\un_1,\ve, V,\Lambda_1)\right| & \le & \frac{\lambda^{-k/2}}{\un_1!} \int_{\real^d} 
\varphi_k(\sqrt{\lambda}\,\ve V\uy_1)\, \left|H_{\un_1}(\uy_1)\right|\, \varphi_d(\uy_1)\,d\uy_1   \\
&\le & \frac{K_1^d\lambda^{-k/2}}{\sqrt{\un_1!}} \int_{\real^d} \varphi_k(\sqrt{\lambda}\,\ve V\uy_1)\, d\uy_1 = \frac{K_1^d\lambda^{-k/2}}{\sqrt{\un_1!}}
\end{eqnarray*}
where we have used the following inequality:\; 
$\displaystyle \sup_x \left| H_l(x)\varphi(x) /\sqrt{l!} \right| \le K_1$, for some constant $K_1$ (see \cite{szego}).

Next, noticing that $\varphi_k(\sqrt{\lambda}\,\ve V\uy_1)$ converges to a Dirac 
delta on $V^{\perp}$ as $\ve\to 0$, we can then write
\begin{equation}\label{eqn:c1-bound}
\sup_{\ve}\,  \,\, c^2_1(\un_1,\ve, V,\Lambda_1) \, \un_1\! !   \,\,\,\, \le \,\,\,\,\,K_1^{2d}\lambda^{-k} .
\end{equation}

\noindent $\underline{\text{{\scshape Computing $c_2$}}}$

The coefficient $c_2(\un_2,u, V, \Lambda_2)$ is the Hermite coefficient of 
$G^V_{2,u}\circ \Lambda_2$ (with $G^V_{2,u}(f_2(x))$ defined in \eqref{eqn:G2-def}), i.e.
\begin{equation}\label{c2}
c_2(\un_2,u, V, \Lambda_2)= \frac1{\un_2 !}\int_{\real^{D-d}} \left(G^V_{2,u}\circ\Lambda_2\right)(\uy_2)\, H_{\un_2}(\uy_2)\,\varphi_{D-d}(\uy_2)\,d\uy_2
\end{equation}
or, equivalently, 
$$
c_2(\un_2,u, V, \Lambda_2) = \frac1{\un_2!} \mathbb{E}\left[ \left(G^V_{2,u}\circ\Lambda_2\right)(\underline{Z}_2)\times H_{\un_2}(\underline{Z}_2)\right]
$$
introducing $\underline{Z}_2$ as a $(D-d)$-dimensional standard normal variable.

Next, using Cauchy-Schwarz inequality, we can conclude that 
\begin{eqnarray*}
c_2(\un_2,u, V, \Lambda_2) 
&\le &  (\un_2!)^{-1/2} \left\{\mathbb{E}\left[\left(G^V_{2,u}\circ\Lambda_2\right)(\underline{Z}_2)\right]^2\right\}^{1/2}\\
& \le &   (\un_2!)^{-1/2} \left(\mathbb{P}[f(x)>u]\right)^{1/2} \, \left\{\mathbb{E}[\det(V\grad^2 f(x)V^T)]^4\right\}^{1/4}.
\end{eqnarray*}
Again using the invariance of $c_2(\un_2,u, V, \Lambda_2)$ with respect to $V$, we can choose $V$ to 
be the line span of $(e_1,\ldots,e_k)$, where $\{e_i\}_{i=1}^d$ is the canonical basis of $\real^d$. 
Writing $\displaystyle \grad^2 f(x)|_{k\times k}$ as the top left $k\times k$ minor of $\grad^2 f(x)$, we have 
$\displaystyle \mathbb{E}[\det(V\grad^2 f(x)V^T)]^4 = \mathbb{E}[\det(\grad^2 f(x)|_{k\times k})]^4$. 
Then, using Wick's formula, we can obtain an upper bound for $\displaystyle \mathbb{E}[\det(V\grad^2 f(x)V^T)]^4$. \\
On the other hand, 
$\mathbb{P}[f(x) >u]$ can be bounded above (and below) by the standard Mill's ratio, implying there exists
$K_{2,u}\in (0,\infty)$ such that
\begin{equation}\label{eqn:c2-bound}
c_2(\un_2,u, V, \Lambda_2)\le K_{2,u}.
\end{equation}
\begin{rk}
Now that we have seen precise expressions for the Hermite \, coefficients $c_1(\un_1, V,\Lambda_1)$ and 
$c_2(\un_2, u, V,\Lambda_2)$, and we understand that these coefficients do not depend on the choice of $V$,
therefore, we shall replace $V$ by its dimension $k$ in the above notations.
In particular, we shall now redefine $$c_1(\un_1, k,\Lambda_1)\definedas c_1(\un_1, V,\Lambda_1),$$
$$
c_2(\un_2, u, k,\Lambda_2)\definedas c_2(\un_2, u, V,\Lambda_2),$$
and $$c(\un, u, k,\Lambda)\definedas c_1(\un_1, k,\Lambda_1) \times c_2(\un_2, u, k,\Lambda_2).
$$
\end{rk}
With these notations, and armed with the fact that 
$\displaystyle \mathbb{E}\left[\left(\phi_k(\partial_{k} \TT_{V^*})\right)^2 \right] < \infty$, we can 
conclude that the following infinite expansion holds in $L^2$
\begin{eqnarray}
\phi_k(\partial_{k} \TT_{V^*}) & = & \sum_{q=0}^{\infty} \underset{\un\in\mathbb{N}^D; |\un|=q}{\sum} 
c(\un,u,k,\Lambda) \int_{\TT_{V^*}} H_{\un}(Z(x))\,dx \quad \label{eqn:phi-k-Hn} \\
&\definedas & \sum_{q=0}^{\infty} \,\, J_q(\phi_k(\partial_{k} \TT_{V^*})),  \label{eqn:phi-k-Hermite}
\end{eqnarray}
where $J_q(\phi_k(\partial_{k} \TT_{V^*}))$ is the projection of $\phi_k(\partial_{k} \TT_{V^*})$ onto the $q$-th chaos.
In addition, we have the following expansion for $\LL_{d-k,k}\left(A_u(f;\TT)\right)$.

\begin{prop}\label{prop:hermite}
For $f$ satisfying the assumptions set forth in Section \ref{subsec:assume}, the following expansion holds in $L^2$:
\begin{equation}\label{eqn:Lk-Hermite}
\LL_{d-k,k}\left(A_u(f;\TT)\right) = \sum_{q=0}^{\infty} \underset{\un\in\mathbb{N}^D; |\un|=q}{\sum} 
c(\un,u,k,\Lambda) \int_{\Gr(d,k)} \int_{\TT_{V^*}} H_{\un}(Z(x))\,dx\,d\lambda^d_k(V^*).
\end{equation}
\end{prop}

{\bf Proof:} First consider the finite sum
\begin{equation}\label{eqn:LL-Q}
\LL^{(Q)}_{d-k,k}\left(A_u(f;\TT)\right) = \sum_{q=0}^{Q} \underset{\un\in\mathbb{N}^D; |\un|=q}{\sum} c(\un,u,k,\Lambda) \int_{\Gr(d,k)}
 \int_{\TT_{V^*}} H_{\un}(Z(x))\,dx\,d\lambda^d_k(V^*).
\end{equation}
Also, define $\phi^{(Q)}_k(\partial_{k} \TT_{V^*})$ as the projection of $\phi_k(\partial_{k} \TT_{V^*})$ onto the first $Q$ orders of the Hermite expansion given in \eqref{eqn:phi-k-Hermite}. Then we can write
$$
\LL^{(Q)}_{d-k,k}\left(A_u(f;\TT)\right) = \int_{\Gr(d,k)} \phi^{(Q)}_k(\partial_{k} \TT_{V^*}) \,d\lambda^d_k(V^*).
$$
Writing $\|\cdot\|_2$ for $L^2$ norm, we have,
\begin{eqnarray*}
&& \left\|\LL_{d-k,k}\left(A_u(f;\TT)\right) - \LL^{(Q)}_{d-k,k}\left(A_u(f;\TT)\right)\right\|^2_2 \\
&= & \left\| \int_{\Gr(d,k)} \left[ \phi^{(Q)}_k(\partial_{k} \TT_{V^*}) - \phi_k(\partial_{k} \TT_{V^*})\right] \,d\lambda^d_k(V^*) \right\|^2_2\\
&\le & \left| \TT\right| \left[ \begin{array}{c} d\\ k\end{array}\right] \int_{\Gr(d,k)} \|\phi^{(Q)}_k(\partial_{k} \TT_{V^*}) - \phi_k(\partial_{k} \TT_{V^*})\|^2_2\, d\lambda^d_k(V^*).
\end{eqnarray*}
Next, using computations similar to those in Section \ref{subsec:sq-int}, we can conclude that there exists a finite, positive $C_{\TT}$ such that $\|\phi_k(\partial_{k} \TT_{V^*})\|^2_2 \le C_{\TT} |\TT_{V^*}|$. 
Notice also, via \eqref{eqn:phi-k-Hermite}, that 
$\displaystyle \|\phi^{(Q)}_k(\partial_{k} \TT_{V^*}) - \phi_k(\partial_{k} \TT_{V^*})\|_2 \underset{Q\to\infty}{\to} 0$. We can then conclude, via the dominated convergence theorem, that
$$\left\|\LL_{d-k,k}\left(A_u(f;\TT)\right) - \LL^{(Q)}_{d-k,k}\left(A_u(f;\TT)\right)\right\|^2_2 \,\,\,\,\underset{Q\to\infty}{\longrightarrow} \,\,\,\, 0$$ 
which proves Proposition \ref{prop:hermite}.  \qed

\begin{rk}\label{rk:orthogonal}~
\begin{itemize}
\item[(i)] Another convenient way of writing $\LL_{d-k,k}(A_u(f;\TT))$ is to express the above expansion as
$$
\LL_{d-k,k}(A_u(f;\TT)) = \sum_{q=0}^{\infty} \int_{\Gr(d,k)}\,\, J_q(\phi_k(\partial_{k} \TT_{V^*}))\,\,d\lambda^d_k(V^*).
$$
\item[(ii)] As a consequence of \cite[Lemma 3.2]{taqqu77}, 
we note that for any $V^*_1,V^*_2 \in\Gr(d,k)$.
$$
\mathbb{E}\left( \int_{\TT_{V_1^*}}\int_{\TT_{V_2^*}} H_{\un}(Z(x))\, H_{\um}(Z(y))\,dx\,dy\right) = 0, \,\,\,\,\,\,\,\text{whenever } |\un| \neq |\um|,
$$
which in turn implies that the expansion in \eqref{eqn:Lk-Hermite} is indeed orthogonal.
\end{itemize}
\end{rk}

\subsection{Variance bounds}
\label{subsec:var-bounds}

Let us define the appropriately normalized quantities of interest
\begin{equation}
\label{eqn:normalized-Lk-1}
\LL^\#_{d-k}(\TT) = \frac{1}{|\TT|^{1/2}}\left( \LL_{d-k}(A_u(f;\TT)) - \mathbb{E}\left[ \LL_{d-k}(A_u(f;\TT))\right]\right) \quad
\end{equation} 
and
\begin{equation}
\label{eqn:normalized-Lk-2}
\LL^\#_{d-k,l}(\TT) = \frac{1}{|\TT|^{1/2}}\left( \LL_{d-k,l}(A_u(f;\TT)) - \mathbb{E}\left[ \LL_{d-k,l}(A_u(f;\TT))\right]\right), \, 0\le l\le k.
\end{equation} 
 
We want to ensure that the variance of $\LL^{\#}_{d-k,k}(A_u(f;\TT))$ converges to a finite
positive quantity as $\TT\to \real^d$, and that the variance of $\LL^{\#}_{d-k,k}(A_u(f;\TT))$ for each $l=0,\ldots,(k-1)$
can be made as small as we wish, by choosing appropriately large set $\TT$.

\begin{prop}\label{prop:var-Lk} ~
With the above notation, the variance of $\LL^\#_{d-k}$ is given by
\begin{equation}\label{eqn:variance-1}
\var\left( \LL^\#_{d-k}\left(\TT\right) \right) = \var\left( \LL^\#_{d-k,k}\left(\TT\right) \right) + {\rm o}(1), \,\,\,\,\,\,\text{ as } \TT\to\real^d. 
\end{equation}
The asymptotic variance of $\displaystyle \LL^\#_{d-k}\left(\TT\right) $, as $\TT\to\real^d$, is finite, 
non zero,  and can be expressed as
\begin{equation}\label{eqn:variance-2}
\lim_{\TT\to\real^d}\var\left(\LL^\#_{d-k,k}\left(\TT\right) \right) = \sum_{q=1}^\infty V_q^k \,\,\, \in (0,\infty)
\end{equation}
where $\displaystyle V_q^k= \lim_{\TT\to {\real}^d} \var\left( \int_{\Gr(d,k)} J_q(\phi^{\#}_k(\partial_k \TT_{V^*}))\,d\lambda^d_k(V^*)\right)$, 
with \\ $\displaystyle \phi^{\#}_k(\partial_k \TT_{V^*})
=\frac1{|\TT|^{1/2}}\,\left(\phi_k(\partial_k \TT_{V^*})-\mathbb{E}\left[ \phi_k(\partial_k \TT_{V^*})\right]\right)$.
\end{prop}

Using the Hermite expansion of $\LL_{d-k,k}(A_u(f;\TT))$ and the orthogonality of the chaos expansion 
(see Remark~\ref{rk:orthogonal} {\it (ii)}), we can formally express the variance of $\LL_{d-k,k}(A_u(f;\TT))$ as
\begin{eqnarray*}
&& \var\left( \LL_{d-k,k}(A_u(f;\TT))\right) \\
&=&\sum_{q=1}^{\infty} \var\left( \underset{\un\in\mathbb{N}^D; |\un|=q}{\sum} c(\un,u,k,\Lambda)\int_{\Gr(d,k)} \int_{\TT_{V^*}} H_{\un}(Z(x))\,dx\,d\lambda^d_k(V^*)\right) \\
&=& \sum_{q=1}^{\infty}\!\!\! \underset{\begin{array}{c}\scriptstyle |\un|=|\um|=q\\ \scriptstyle\un,\um\in\mathbb{N}^D \end{array}}{\sum}
\!\!\!\!\!\!\!\!\! c(\un,u,k,\Lambda)\,c(\um,u,k,\Lambda) \!\!\!\! \!\!\!\!\!\!\! \underset{U^*,V^*\in Gr(d,k)}{\iint} d\lambda^d_k(U^*)\,d\lambda^d_k(V^*) \\
&   & \hskip 1.5cm \times  \left(\underset{\TT_{U^*}}{\int}\underset{\TT_{V^*}}{\int} \mathbb{E}\left[H_{\un}(Z(x))H_{\um}(Z(y))\right]  dx\,dy\right) 
\end{eqnarray*}

The sketch and main arguments (e.g. Arcones bound) to prove Proposition \ref{prop:var-Lk} are given in \cite{KL-jotp}, 
with an extra step for the term $o(1)$ which follows from \cite{EL13}. The main difficulty relies then, once again, in 
the fact that we do not integrate simply on a $d$-dimensional box, but on Grassmanians, which requires tricks to 
circumvent the difficulty of computations. 

{\bf Proof of Proposition~\ref{prop:var-Lk}.} ~
First let us show that $\displaystyle \var\left( \LL^\#_{d-k,k}\left(\TT\right) \right)<\infty$. We have,
using \eqref{eqn:def-LL(d,d-k)}, then \eqref{eqn:phi-k-Hermite},

\begin{eqnarray} \label{Anmu}
&& \var\left( \LL^\#_{d-k,k}\left(\TT\right) \right) \nonumber\\
&=&\frac1{ \left| \TT\right|} \int_{\Gr(d,k)}\int_{\Gr(d,k)} \cov\left( \phi_k(\partial_k\TT_{V^*}),\phi_k(\partial_k\TT_{U^*})\right)\,d\lambda^d_k(V^*)\,d\lambda^d_k(U^*) \nonumber \\
&=& \frac1{ \left| \TT\right|}  \sum_{q=1}^{\infty}\sum_{|\un|=q}\sum_{|\um|=q}  
c(\un,u, k,\Lambda) \,\, c(\um, u, k,\Lambda)\, A(\un,\um,u,k,\TT) 
\end{eqnarray}
with
\begin{eqnarray*}
&&A(\un,\um,u,k,\TT) : = \underset{\Gr(d,k)}{\int} \underset{\Gr(d,k)}{\int} \,d\l^d_k(U^*)\,d\l^d_k(V^*)\\
&& \hskip 1.0cm \times 
\left[\underset{x\in\partial_k \TT_{U^*}}{\int}\underset{y\in\partial_k \TT_{V^*}}{\int}
\mathbb{E}\left[H_{\un}\left(Z(x)\right)H_{\um}\left(Z(y)\right)\right] \, dx \,dy\right].
\end{eqnarray*}
Notice that since $Z(x)$ and $Z(y)$, individually, are standard Gaussian vectors, then using Mehler's formula 
(or equivalently, the diagram formula), we have that for $|\un|=|\um|=q$,
\begin{eqnarray*} 
&& A(\un,\um,u,k,\TT) \\
&=& \sum_{\scriptsize\begin{array}{c} d_{ij}\ge 0 \\ \sum_{i}d_{ij} = n_j \\ \sum_{j}d_{ij} = m_i\end{array}}
\un! \, \um ! \underset{U^*\in\Gr(d,k)} {\int}  \underset{x\in\partial_k \TT_{U^*}}{\int} \,dx \,d\l^d_k(U^*) \\
&& \times \left[\underset{V^*\in\Gr(d,k)}{\int} 
\underset{w\in\partial_k (x-\TT_{V^*})}{\int} \prod_{1\le i,j\le D} \frac{\gamma^{d_{ij}}_{ij}(w)}{d_{ij}!} \,dw \,d\l^d_k(V^*)\right] \\
&\le &  \!\!\!\!\!\!\! \sum_{\scriptsize\begin{array}{c} d_{ij}\ge 0 \\ \sum_{i}d_{ij} = n_j \\ \sum_{j}d_{ij} = m_i\end{array}} \!\!\!\!\!\!\!\!\un! \, \um ! \left[
 \underset{W^*\in\Gr(d,k)}{\int} \underset{w\in\partial_k (2\TT_{W^*})}{\int} \prod_{1\le i,j\le D} \frac{ |\gamma^{d_{ij}}_{ij}(w) | }{d_{ij}!} \,dw \,d\l^d_k(W^*)\right] \\
&& \hskip 1.0cm \times \underset{U^*\in\Gr(d,k)}{\int}\underset{x\in\partial_k (\TT_{U^*})}{\int} \!\!\! \! \,dx \,d\l^d_k(U^*)
\end{eqnarray*}
where the inequality is a result of the observation that for any $x\in\TT$, we have the inclusion 
$(x-\TT_{V^*}) \subset  2\TT_{W^*}$, for some $W^*\in \Gr(d,k)$.
Before proceeding any further, we may observe the following.
\begin{lem}\label{lem:gammaGrk} 
Let $\theta$ be a nonnegative real valued, integrable function 
defined on $\real^d$, then
$$
\underset{V^*\in\Gr(d,k)} {\int}\underset{\; \; w\in\partial_k (2\TT_V^*)}{\int} \!\!\theta(w) \,\,dw \, d\l^d_k(V^*) \le 
\left[ \begin{array}{c} d\\k\end{array}\right] \int_{\real^d} \theta(z) dz.
$$
\end{lem}

{\bf Proof:}
The double integral in question can be bounded from above by first replacing the integral over $\partial_k(2\TT_{V^*})$ 
by integral over $V^*$. Then, since
$\Gr(d,k)$ is isometric to $\text{Gr}(d,k)\times \real^{d-k}$, we note that, for any fixed $V\in\text{Gr}(d,k)$, we have
$\cup_{x\in\real^{d-k}}(V+x) = \real^d$. Therefore the above integral can be bounded above by
$$\int_{V\in\text{Gr}(d,k)} \left(\int_{\real^d} \theta(w)\,dw \right) \, d\sigma^d_k(V)$$
where $\sigma^d_k$ is the invariant measure on the Grassmannian $\text{Gr}(d,k)$ such that
$$\sigma^d_k\left( \text{Gr}(d,k)\right) = \left[ \begin{array}{c} d\\k\end{array}\right],$$
which proves the assertion of the lemma.
\qed

In view of Lemma \ref{lem:gammaGrk}, an upper bound for $A(\un,\um,u,k,\TT)$ can be obtained as
\begin{eqnarray*}
&& A(\un,\um,u,k,\TT) \\
&\le & \left[ \begin{array}{c} d\\k\end{array}\right]^2
\!\!\!\!\!\!\!\!\!\!\!\!\!\sum_{\scriptsize\begin{array}{c} d_{ij}\ge 0 \\ \sum_{i}d_{ij} = n_j \\ 
\sum_{j}d_{ij} = m_i\end{array}}\!\!\!\!\!\! \un!\, \um!
\left[
 \underset{w\in \real^d}{\int} \prod_{1\le i,j\le D} \frac{\left|\gamma^{d_{ij}}_{ij}(w)\right|}{d_{ij}!} dw \right]
\underset{x\in \TT}{\int} \,dx
\end{eqnarray*}
where in the second integral we have used the Crofton formula \eqref{crofton-volume}.

Further, under hypothesis (H2), and for $|n|=|m|=q$, there exists a constant $C^*$ such that 
$$
\sum_{\scriptsize\begin{array}{c} d_{ij}\ge 0 \\ \sum_{i}d_{ij} = n_j \\ \sum_{j}d_{ij} = m_i\end{array}}\!\!\!\!\!\! \un!\um! \prod_{1\le i,j\le D} \! \frac{\left|\gamma^{d_{ij}}_{ij}(w)\right|}{d_{ij}!} \, \le \, C^*\,\psi^q(w).
$$

Therefore we obtain that,  for $|n|=|m|=q$, 
$$
A(\un,\um,u,k,\TT) \le  C^* |\TT|\,\left[ \begin{array}{c} d\\k\end{array}\right]^2  \int_{\real^d}  \psi^q(w) \,dw\, .
$$

Next, we prove, asuming Arcones condition (H2), 
that $|\TT|^{-1}A(\un,\um,u,k,\TT)$ converges as $\TT\to\real^d$. We shall check that it is
Cauchy in the parameter $T$, the edge length of $\TT$. 
Let us take boxes $\TT_1 = [-T_1,T_1]^d$ and $\TT_{1,2} = [-(T_1+T_2),T_1+T_2]^d$, and prove
that 
$$\left|\frac{A(\un,\um,u,k,\TT_1)}{|\TT_1|} - \frac{A(\un,\um,u,k,\TT_2)}{|\TT_2|}\right| \to 0 \,\,\,\,\,\,\text{ as }
\,\,\,\, T_1,T_2\to\infty$$
Clearly, 
\begin{eqnarray*}
&& \left|\frac{A(\un,\um,u,k,\TT_{1,2})}{|\TT_{1,2}|} - \frac{A(\un,\um,u,k,\TT_1)}{|\TT_1|}\right| \\
&\le & \frac1{|\TT_{1,2}|}\left| A(\un,\um,u,k,\TT_{1,2}) - A(\un,\um,u,k,\TT_1)\right| + 
\frac{A(\un,\um,u,k,\TT_1)}{|\TT_1|}\left|\frac{|\TT_1|}{|\TT_{1,2}|} - 1\right|\\
&=& \frac1{|\TT_{1,2}|}\left| A(\un,\um,u,k,\TT_{1,2}) - A(\un,\um,u,k,\TT_1)\right| + 
\frac{A(\un,\um,u,k,\TT_1)}{|\TT_1|}\left|\frac{T^d_1}{(T_1+T_2)^d} - 1\right|\\
&& :=  I + II
\end{eqnarray*}
Clearly, the coefficient in $II$ can be bounded uniformly as a result of previous computations, and the volume terms 
converge to zero as $T_1$ increases to infinity. For part $I$, notice that the difference
$\left|A(\un,\um,u,k,\TT_{1,2}) - A(\un,\um,u,k,\TT_1)\right|$ can be shown to be of the same order as
$$
\left|\TT_{1,2} \backslash \TT_1\right| \int_{\TT_{1,2} \backslash \TT_1} \psi^q(w)dw\,.
$$
The coefficient of the integral above, when compared with $|\TT|_{1,2}^{-1}$, converges to one.
However, since the domain of integration escapes to infinity, the integral converges to zero due to 
integrability of $\psi^q$.

Hence, we can conclude that the sequence $|\TT|^{-1}A(\un,\um,u,k,\TT)$ is 
Cauchy in the variable $\TT$, meaning that, for $|\un|=|\um|$,
$$
|\TT|^{-1}A(\un,\um,u,k,\TT) \to A(\un,\um,u,k) \,\,\, \text{ as } \TT\to \real^d \,\,\,\,(\text{or, equivalently, as } T\to\infty)
$$
where the limit $A(\un,\um,u,k)$, using the arguments of Lemma \ref{lem:gammaGrk}, can be identified as
\begin{equation}\label{eqn:Anmu}
A(\un,\um,u,k) = \left[ \begin{array}{c}d \\ k\end{array}\right]^2 \un! \, \um! \!\!\!\!
\sum_{\scriptsize\begin{array}{c} d_{ij}\ge 0 \\ \sum_{i}d_{ij} = n_j \\ \sum_{j}d_{ij} = m_i\end{array}}
\int_{\real^d} \underset{1\le i,j\le D}{\Pi} \frac{\gamma^{d_{ij}}_{ij}(w)}{d_{ij}!}\,
\,dw,
\end{equation}
which, in turn implies that 
$$
\var\left( \int_{\Gr(d,k)} J_q \left( \phi^{\#}_k(\partial_k \TT_{V^*})\right)d\l^d_k(V^*)\right) \rightarrow V^k_q, \quad\text{as}\quad\TT\to\real^d.
$$
\hskip 0.5cm

\noindent$\underline{\text{{\scshape Finiteness of the limiting variance}}}$

We shall proceed as usual (see \cite{KL-SPA} or \cite{KL-jotp}). Introducing  $\Pi_Q\left( \LL^{\#}_{d,d-k}(A_u(f;\TT))\right)$ 
as the projection of $\LL^{\#}_{d,d-k}(A_u(f;\TT))$ onto the first $Q$ chaos, we shall show that
\begin{equation}\label{eqn:var-unif}
\var\left( \LL^{\#}_{d,d-k}(A_u(f;\TT)) - \Pi_Q\left( \LL^{\#}_{d,d-k}(A_u(f;\TT))\right)\right) 
\underset{Q\to\infty}{\longrightarrow} 0, \text{ uniformly in } \TT,
\end{equation} 
and conclude the finiteness of the limiting variance by a simple application of Fatou's lemma.

Let us begin with observing that $\LL_{d,d-k}(A_u(f;\TT))$ is an additive set functional. In particular, the set $\TT$ can
be written, as in \cite{EL13} as a union of disjoint unit cuboids (w.l.o.g. let $T$ be integer). Therefore,
$\LL_{d,d-k}(A_u(f;\TT))$ can be written as a sum of a stationary sequence of random variables
 where these random variables are an evaluation of $\LL_{d,d-k}(A_u(f;\cdot)$ on $[0,1)^d$, and its various integer shifts.

Next invoking stationarity of the field $(\grad f, \grad^2 f, f)$ (and $Z$), we know that the variance of the sum of a stationary sequence
is of the order of the cardinality of the sum if the covariance decays at an appropriate rate. Using this precise argument,
and following the computations of \cite{EL13}, we can conclude \eqref{eqn:var-unif}. In following the arguments of
\cite{EL13}, it is important to note that our estimates for the coefficients in the Hermite expansion match with those in \cite{EL13}.\\

Now we shall show that the variance corresponding to lower dimensional faces of $\TT_{V^*}$, is indeed $o(1)$
for large $\TT$ as expressed in Proposition \ref{prop:var-Lk}.

Recall the decomposition of $\LL_{d-k,k}$ from equation \eqref{eqn:L0-sum-dim}. Then,
\begin{eqnarray*}
&& \var\left( \LL_{d-k}(A_u(f;\TT))\right)\\
 &=&\sum_{l=0}^k \var\left( \LL_{d-k,l}(A_u(f;\TT))\right)
+ 2\sum_{l<m} \cov\left( \LL_{d-k,l}(A_u(f;\TT)),\LL_{d-k,m}(A_u(f;\TT))\right).
\end{eqnarray*}
It suffices to show that $\var\left( \LL_{d-k,l}(A_u(f;\TT))\right) = o(|\TT|)$ for each $l=0,\ldots,(k-1),$
in order to conclude the second part of the assertion in equation \eqref{eqn:variance-1}.

Let us define
\begin{equation}\label{eqn:Rdk}
R(d,k,l,\TT) = \frac1{\left| \TT\right|^{1/2}} \underset{\text{Graff}(d,k)}{\int} |\partial_l \TT_{V^*}|^{1/2} \phi^\#_l\left( \TT_ {V^*}\right) \,d\lambda^d_k(V^*).
\end{equation}
In view of $\partial_l \TT = \!\!\!\!\underset{V^*\in\Gr(d,k)}{\bigcup} \partial_l \TT_{V^*}$, and the above computations leading to 
$A(\un,\um,u,k)$, we note that $\var(R(d,k,l,\TT)) $ can be shown to be $O(|\partial_l \TT|)$, or equivalently $O(T^l)$
under the assumption (H2), implying that the lower dimensional faces, asymptotically, do not contribute to the variance
of $\LL^{\#}_{d-k}(A_u(f;\TT))$.

\hskip 0.5cm 

\noindent$\underline{\text{{\scshape Nondegeneracy of the limit}}}$

Finally, it remains to show that $\displaystyle \lim_{\TT\to\real^d}\var\left( \LL^{\#}_{d-k,k}(A_u(f;\TT))\right)>0$.
Using the orthogonality of chaos, it suffices to show that $V^k_1 >  0$.

First, we shall simplify the expression for $V^k_1(\TT)$ by introducing the canonical basis $(\ue_i)_{1\le i\le D}$ of 
$\real^D$ in \eqref{Anmu}, and writing
$$
V^k_1(\TT) =\frac1{|\TT|} \sum_{i=1}^D \sum_{j=1}^D c(\ue_i,u, k,\Lambda) \,c(\ue_j, u, k,\Lambda)\, A(\ue_i,\ue_j,u, k,\TT). 
$$
Then, writing $\{\ue_{i1}\}_{i=1}^d$, $\{\ue_{j2}\}_{j=1}^{D-d}$ for canonical basis of dimension $d$, $(D-d)$ respectively, and observing that 
$c_1(\ue_{i1},k,\Lambda_1) = 0$ (by \eqref{eqn:c1-precise}), 
the limiting variance corresponding to the first chaos, again using equation \eqref{eqn:c1-precise} 
for precise expression of $c_1(\underline{0},k,\Lambda_1)$, is given by
\begin{equation}\label{eqn:V_1-lower}
V^k_1 = (2\pi\lambda)^{-k}\sum_{i=d+1}^D \sum_{j=d+1}^D \, c_2(\ue_{i2},u,k,\Lambda)\, c_2(\ue_{j2},u,k,\Lambda) \,A(\ue_{i2},\ue_{j2},u,k) 
\end{equation}
with  $A(\ue_{i2},\ue_{j2},u,k)$ as defined in \eqref{eqn:Anmu}, given by
$$
A(\ue_{i2},\ue_{j2},u,k) = \left[ \begin{array}{c}d \\ k\end{array}\right]^2 
\int_{\real^d} \gamma_{\ue_{i2},\ue_{j2}}(w)\, dw,
$$
where $\gamma_{\ue_{i2},\ue_{j2}}$ denotes the covariance function corresponding to the pair of indices 
which correspond to the position of $1$'s in $(\underline{0},\ue_{i2})$ and $(\underline{0},\ue_{j2})$, respectively,
where $\underline{0}$ is a $d$-dimensional row vector of zeros.

As in \cite{EL13}, we have that $A(\ue_{i2},\ue_{j2},u,k)=0$, whenever $(i,j)\neq(D,D)$, thus further simplifying \eqref{eqn:V_1-lower} to
\begin{equation}\label{eqn:VK1T}
V^k_1 =   c^2_2(\ue_{D2},u,k,\Lambda_2) \,A(\ue_{D2},\ue_{D2},u, k). 
\end{equation}

We shall estimate separately the two terms appearing above.

Let us begin with $c_2(\ue_{D2},u,k,\Lambda_2)$ given in \eqref{c2}. We have
\begin{eqnarray}\label{eq:c2V}
&& c_2(\ue_{D2},u,k,\Lambda_2) \nonumber\\
& = & \int_{\real^{D-d}} \left(G_{2,u}^V\circ \Lambda_2\right) (\uy_2) H_{\ue_{D2}}(\uy_2) \varphi_{D-d}(\uy_2) d\uy_2  \nonumber\\
&=& \int_{\real^{D-d-1}} \varphi_{D-d-1}(\uy^*_2) \left( \int_{\real} \left(G_{2,u}^V\circ \Lambda_2\right) (\uy_2) \, y_{2D}\, \varphi(y_{2D})dy_{2D}\right) d\uy^*_2.\,\qquad
\end{eqnarray}
Let us consider the lower triangular matrix $\Lambda_2$ such that its first element $(\Lambda_2)_{11}$ equals $1$ 
(as in \cite{EL13}),  {\it i.e.} of the form 
$\displaystyle \Lambda_2= \left(\begin{array}{cc} L &  \bf{0}\\ \ug^T & l \end{array} \right),$ 
with $L$ a lower triangular $(D-d-1)\times(D-d-1)$ matrix, $\ug^T$  a $1\times (D-d-1)$ matrix, and $l>0$.
With the above notation, we can write
\begin{equation}\label{eq:G2det}
G_{2,u}^V\circ \Lambda_2 (\uy_2) = \det(V\mathcal{M}(L\uy^*_2)V^T)\, 1_{\left\{\ug\,\uy^*_2 + l \,y_{2D}\ge u\right\}}
\end{equation}
where $\mathcal{M}(L\uy^*_2)$ is the symmetric matrix obtained by appropriately arranging the elements of the vector 
$L\uy^*_2$.

We can certainly think of the map $\uy^*_2 \mapsto \mathcal{M}(L\uy^*_2)$ as a linear map, therefore, there exists $a^b_{ij}$ such that
$$
\mathcal{M}(L\uy^*_2) _{ij}= \sum_{b=1}^{D-d-1} a^b_{ij} \,\uy^*_{2,b}.
$$
Again recalling that $c_2$ does not depend on the choice of $V$, we shall fix the matrix $V$ as 
$\left[ I_k \,;\, \underline{0}\right]$, where $I_k$ is $k\times k$ identity matrix and $\underline{0}$ is a 
$k\times (d-k)$ matrix of zeros. Then,
\begin{equation}\label{eq:keyM}
\left(V\mathcal{M}(L\uy^*_2) V^T\right)  = \mathcal{M}(L\uy^*_2)\vert_{k\times k}
\end{equation}
where the right side is the notation for the top left $k\times k$ minor of $\mathcal{M}(L\uy^*_2)$.

This latter argument \eqref{eq:keyM} is key, since the next computations will then be similar as those done in a $d$-dimensional box (\cite{EL13}). We now give a brief sketch of the major steps involved to provide an overview of the full computation. We have
$$
\det\left(\mathcal{M}(L\uy^*_2)|_{k\times k}\right) = \sum_{\sigma\in \mathcal{S}_k} \text{sgn}(\sigma) \prod_{i=1}^k 
\left[ \sum_{b=1}^{D-d-1} a^b_{i\sigma(i)} \,\uy^*_{2,b}\right].
$$
Subsequently, using arguments similar to those in Lemma A.2 of \cite{EL13} together with isotropy, we can obtain a Hermite
expansion for the determinant as follows
\begin{equation}\label{eqn:hermite-det}
\det(V\mathcal{M}(L\uy^*_2)V^T) = \underset{\un\in\mathbb{N}^{D-d-1}: |\un|=k}{\sum} \alpha_{\un}(L,V) H_{\un}(\uy^*_2).
\end{equation}
Combining \eqref{eq:c2V}, \eqref{eq:G2det} and \eqref{eqn:hermite-det}, and using that $\displaystyle y\varphi(y)=-\varphi'(y)$ to compute the integrand on $y_{2D}$, we obtain the following (for more details, we refer the reader to the proof of Lemma 2.2 of \cite{EL13})
\begin{eqnarray*}
&& c_2(\ue_{D2},u,k,\Lambda_2) \\
&=&\underset{\un\in\mathbb{N}^{D-d-1}: |\un|=k}{\sum} \alpha_{\un}(L,V) 
\int_{\real^{D-d-1}} H_{\un}(\uy^*_2) \varphi\left(\frac1l \left(u - \langle\ug,\uy^*_2 \rangle \right)\right) \varphi_{D-d-1}(\uy^*_2)\,d\uy^*_2\\
&=& \underset{\un\in\mathbb{N}^{D-d-1}: |\un|=k}{\sum} \alpha_{\un}(L,V) (-1)^k 
\int_{\real^{D-d-1}}  \varphi\left(\frac1l \left(u- \langle \ug,\uy^*_2 \rangle \right)\right) \varphi_{D-d-1}^{(n)}(\uy^*_2)\,d\uy^*_2 \\
&=& l\,H_k(u)\,\varphi(u) \underset{\un\in\mathbb{N}^{D-d-1}: |\un|=k}{\sum} \alpha_{\un}(L,V)H_{\un}(\ug) \\
&=& l\,H_k(u)\varphi(u)\,\det\left(V\mathcal{M}(L\ug) V^T\right).
\end{eqnarray*}

Since we can write $\displaystyle \mathcal{M}(L\ug)= -\lambda I_k$ with $\lambda=-r_{ii}(0)$, then
\begin{equation}\label{eq:claim1}
 c_2(\ue_{D2},u,k,\Lambda_2) = l\,H_k(u)\varphi(u) (-\lambda)^k 
\end{equation}
by way of choosing $V = [ I_k \,;\, \underline{0}]$. \\

Moreover, as in \cite{EL13}, we can write
\begin{equation}\label{eq:claim2}
\underset{w\in \real^D}{\int} \gamma_{\ue_{D2},\ue_{D2}}(w) \,dw = (2\pi)^d\, h(0)\,l^{-2}
\end{equation}
where we recall that $h(0)$ is the spectral density of the field $f$ evaluated at $0$.

Finally, putting together the estimates obtained in \eqref{eq:claim1} and \eqref{eq:claim2} in the following
\begin{eqnarray*}
V_1^k
&=&  (2\pi)^{-k}\,l^2 H_k^2(u)\varphi^2(u)\lambda^{k} \left[ \begin{array}{c}d\\ k\end{array}\right]^2 \underset{w\in \real^D}{\int} \gamma_{\ue_{D2},\ue_{D2}}(w) \,dw,
\end{eqnarray*}
we obtain
\begin{eqnarray*}
V_1^k &=& (2\pi)^{d-k}H_k^2(u)\varphi^2(u)\lambda^{k} h(0)\, \left[ \begin{array}{c}d\\ k\end{array}\right]^2
\end{eqnarray*}
from which we deduce that \, $\displaystyle  V^k_1 \, >\, 0$, hence the second part of Proposition \ref{prop:var-Lk}.
\qed

\subsection{Extension of Breuer-Major theorem to affine Grassmannian case}
\label{subsec:BM}
 
Here we just give a sketchy recall of the literature on CLTs of Breuer-Major type, that can be found in 
\cite{noupebook}, \cite{NPP}. 

In 1983, Breuer-Major provided a CLT for a 1-dimensional centered stationary Gaussian sequence indexed by
$\mathbb{Z}^{\nu}$ for $\nu\ge 1$, satisfying some condition 
on its correlation function. This result was first extended by Giraitis and Surgailis \cite{gir-sur} when considering a continuous 
time setting, then by Arcones \cite{Arcones} with a powerful result holding for vector valued random sequences. 
The proof,  in the discrete case, is 
based on the method of cumulants and diagram formulae. Estrade and L\'eon rewrote it explicitly (see \cite{EL13}, Proposition 2.4) 
in the continuous case following the Nourdin et al's proof (\cite{NPP}) based on Malliavin calculus. To avoid mimicking 
the proof a second time, we shall point out the main quantities that deserve some care, 
due to our general setting.  Let us first state Breuer-Major theorem in this setting.
\begin{prop}\label{prop:BM} 
Let $\TT$ be a $d$-dimensional box $[-T,T]^d$, and let $f$ be a mean zero, unit variance, isotropic Gaussian random field defined 
on $\real^d$ with $C^3$ trajectories. Under the assumptions (H1) to (H3), for any positive integer $Q$, the projection onto the first 
$Q$ chaos $\Pi^Q\left(\LL_{d-k,k}^\#(A_u(f;\TT))\right)$ satisfies
$$
\Pi^Q\left(\LL_{d-k,k}^\#(A_u(f;\TT))\right) \;\overset{d}{\longrightarrow} \; \NN\left(0,\sum_{q=1}^Q V^k_q\right) \quad \text{as} \; \TT\to\real^d,
$$
where  $\displaystyle V_q^k$ is defined in Proposition~\ref{prop:var-Lk}.
\end{prop}

Indeed, we have
\begin{eqnarray*}
&&\frac1{\left| \TT\right|^{1/2}} \underset{\text{Graff}(d,k)}{\int} \phi_k\left( \partial_k\TT_ {V^*}\right) \,d\lambda^d_k(V^*) 
= \frac1{\left| \TT\right|^{1/2}} \underset{\text{Graff}(d,k)}{\int} \underset{\partial_k\TT_V}{\int}  G^V (f_1,f_2)(x)\,dx \,d\lambda^d_k({V^*}) \\
&&= \frac1{\left| \underset{\text{Graff}(d,k)}{\int} \underset{\partial_k\TT_{V^*}}{\int} dx\, d\lambda^d_k({V^*})\right|^{1/2}} \underset{\text{Graff}(d,k)}{\int} \underset{\partial_k\TT_{V^*}}{\int}  G^V (f_1,f_2)(x)\,dx \,d\lambda^d_k({V^*}).
\end{eqnarray*}
where $f_1=\grad f$, $f_2 = (\grad^2 f,f)$, and $G^V(f_1,f_2) = G^V_1(f_1) \times G^V_{2,u}(f_2)$ 
as defined in \eqref{eqn:G1-def} and \eqref{eqn:G2-def}.

Considering the projection onto the first $Q$ chaos, $\Pi^Q\left(\LL_{d-k,k}^\#(A_u(f;\TT))\right)$, defined in \eqref{eqn:var-unif}, we can write, as in the proof of Theorem 2.2 in \cite{NPP} (or in \cite{EL13}), 
\begin{equation}\label{eqn:Wiener-Ito}
\Pi^Q\left(\LL_{d-k,k}^\#(\TT)\right) =\sum_{q=1}^Q I_q(g_{k,q}^{\TT})
\end{equation}
where  $I_q(f)$ denotes the multiple Wiener-It\^o integral (of order q) of $f$ with respect to $W$, and
\begin{eqnarray}\label{g-qk}
g_{k,q}^{\TT} &=& \frac{b^k_{\um}}{\left| \underset{\text{Graff}(d,k)}{\int} \underset{\partial_kT_{V^*}}{\int} dx\, d\lambda^d_k(V^*)\right|^{1/2}}\nonumber\\
&& \hskip 0.2cm \times \underset{\text{Graff}(d,k)}{\int} \underset{\partial_k\TT_{V^*}}{\int}\sum_{\um\in \{1, 2,\cdots,D\}^q}  \, u_{x,m_1}\,\otimes\cdots \,\otimes\, u_{x,m_q}\,dx \, d\lambda^d_k(V^*) \qquad
\end{eqnarray}
 where $b^k_{\um}$ are such that the mapping $\um \rightarrow b^k_{\um}$ is symmetric on $\{1,\cdots, D\}^q$,
 and we have again used isotropy to observe that $b^k_{\um}$ depends on $V^*$ only through its dimension, which is $k$.
Moreover, the functions $(u_{x,j})_{1\le j\le D}$ are orthogonal in $L^2(\real^d)$ such that for the field $Z(x)$ defined in \eqref{Z},
$$Z_j(x) = \int_{\real^d} u_{x,j}(w)\,dW(w)$$
where $W$ is the complex Brownian measure on $\real^d$.

Note that, in writing \eqref{eqn:Wiener-Ito}, we have used the Fubini theorem to interchange the Wiener-It\^o integral and
the integral over the space $\displaystyle \underset{V^*\in\Gr(d,k)}{\cup} \partial_k\TT_{V^*}$.

In order to prove the CLT of $\displaystyle \Pi^Q\left(\LL_{d-k,k}^\#(A_u(f;\TT))\right)$, it is enough to check that, for 
$1\le p,q\le Q$, (see \cite{noupebook} or \cite{NPP}, and for the notation, \cite{EL13})
$$
|| \delta_{pq} V_p^k -\frac1q \left\langle \textsc{D}I_p(g_{k,p}^{\TT}),\textsc{D}I_q(g_{k,q}^{\TT} )\right\rangle_\HH||_2 \;\rightarrow\; 0 
\; \; \text{as}\; \TT\to {\real}^d ,
$$
where $\displaystyle V_q^k$ is defined in Proposition~\ref{prop:var-Lk}, 
and $\textsc{D}$ denotes the Malliavin derivative.
Standard analysis as in \cite{NPP} can be invoked to conclude that it suffices to check that, for $p\le q$, 
$$
|| \frac1q \left\langle\textsc{D}I_p(g^{\TT}_{k,p}),\textsc{D}I_q(g^{\TT}_{k,q} )\right\rangle_\HH||_2 \;\rightarrow\; 0 \; \; \text{as}\; \TT\to {\real}^d  ,
$$
which holds since, on one hand, for the case $p=q$ we have $\displaystyle  || g^{\TT}_{k,q}||^2_{{\HH}^q} = V_q^k(\TT)$ 
which is shown to converge to $V^k_q$ in Proposition~\ref{prop:var-Lk}.
On the other hand, the $e$-th contraction of $g^{\TT}_{k,p}$ satisfies, for $e<p$, 
$$
|| g^{\TT}_{k,p} \, \underset{e}{\otimes}\, g^{\TT}_{k,p}||^2_{{\HH}^{2(p-e)}} \le 
\left( C^p \sum_{\um\in \{1, 2,\cdots,D\}^p}  |b_{\um} |^2\right)^2 \Psi (k)
$$
with some constant $C$, and under (H2), 
\begin{eqnarray*}
&&\Psi(k):=\\
&& \!\!\! \underset{(\text{Graff}(d,k))^4}{\int} \underset{\partial_k\TT_{V^*_1}}{\int} \!\!\!\!\!\cdots\!\!\!\!\!\underset{\partial_k\TT_{V^*_4}}{\int} \!\!\!
\frac{\psi^{e}(t_1-t_2)\psi^{e}(t_3-t_4)\psi^{p-e}(t_1-t_3)\psi^{p-e}(t_2-t_4)}{\left| \underset{\text{Graff}(d,k)}{\int} \underset{\partial_k\TT_{V^*}}{\int} dx\, d\lambda^d_k(V^*)\right|^2}
\,\, \prod_{i=1}^4 dt_i\,d\lambda^d_k(V^*_i).
\end{eqnarray*}
As in \cite{EL13}, we note that $\psi^{e}(t_3-t_4)\psi^{p-e}(t_1-t_3) \le \psi^{p}(t_3-t_4)+\psi^{p}(t_1-t_3)$,
and by Lemma \ref{lem:gammaGrk}, we have 
$$\underset{\text{Graff}(d,k)}{\int} \underset{\partial_k\TT_{V^*}}{\int} \psi^{p}(t_1-t_3)\,ds_3 < \left[ \begin{array}{c} d\\ k \end{array}\right] \int_{\real^d} \psi(z)\,dz < \infty$$
which matches the estimates of \cite{EL13}, and thus we can follow the rest of the arguments verbatim to conclude that for some finite, combinatorial constant $C(k)$, we have
\begin{eqnarray*}
\Psi(k) &\le&
C(k)\, \frac{\left| \underset{\text{Graff}(d,k)}{\int} \underset{\partial_k\TT_{V^*}}{\int} dx\, d\lambda^d_k(V^*)\right|}
{\left| \underset{\text{Graff}(d,k)}{\int} \underset{\partial_k\TT_{V^*}}{\int} dx\, d\lambda^d_k(V^*)\right|^2}
= \, \frac{C(k)}{\left| \underset{\text{Graff}(d,k)}{\int} \underset{\partial_k\TT_{V^*}}{\int} dx\, d\lambda^d_k(V^*)\right|} \; \to \, 0\;\,\,\,
\mbox{ as }\, 
\TT\to {\real}^d.
\end{eqnarray*}
This concludes the proof of Proposition \ref{prop:BM}.
\qed

Collating Propositions  \ref{prop:var-Lk} and \ref{prop:BM} leads to the main result, that is the following CLT 
$$\frac{\LL_{d-k}(A_u(f;\TT))- \mathbb{E}\left[ \LL_{d-k}(A_u(f;\TT))\right]}{|\TT|} \,\,\,\underset{\TT\to\real^d}{\longrightarrow}
 \,\,\,N(0,\sigma^2_{d-k}(u)),$$
where $\sigma^2_{d-k}(u)$ is given by $\sum_{q\ge 1} V^k_q$ in \eqref{eqn:variance-2}.

\begin{rk}
The assumption of isotropy was crucial to circumvent a direct computation of Hermite coefficients
of the LKCs, providing bounds independent of the choice of $V$. Nevertheless the CLT should hold true
under the assumption of stationarity together with hypotheses $(H1)$, $(H2)$ and $(H3)$.
\end{rk}

\section{Discussion} \label{sec:discussion}

{\scshape Extension to general parameter spaces:}

Notice that the only place where we required the box type shape of the parameter space is when we get an upper bound on the limiting variance of $\LL^{\#}_{d-k,k}(A_u(f;\TT))$. 
However, this can be overcome by a limiting procedure. 

Let us partition the space $\real^d$ into small cuboids of volume $\eta$. We can identify these small cuboids by the 
centre of the cuboids. Let $C^{\eta}_T$ be the set of cuboids which completely lie in the set $\TT$, and $B^{\eta}_{\TT}$ 
be the cuboids which have non empty intersection with the set $\TT$ and the complement of $\TT$.

Denoting $\PP_{i,\eta}$ for the elements of the partition of $\real^d$ into cuboids of volume $\eta$, we have
 \begin{eqnarray*}
 && \LL^{\#}_{d,d-k}(A_u(f;\TT)) \\
 &=& \sum_{\PP_{i,\eta}\in C^{\eta}_{\TT}} \LL^{\#}_{d,d-k}(A_u(f;\PP_{i,\eta}))
 + \sum_{\PP_{i,\eta}\in B^{\eta}_{\TT}} \LL^{\#}_{d,d-k}(A_u(f;\PP_{i,\eta})\cap \TT) \\
 &\definedas& \,\, \LL^{\#}_{d,d-k}(A_u(f;\TT),1) + \LL^{\#}_{d,d-k}(A_u(f;\TT),2).
 \end{eqnarray*}
Notice that using stationarity and the decay of covariance function $\gamma$, as in \cite{EL13}, we can conclude that
$$\var\left( \LL^{\#}_{d,d-k}(A_u(f;\TT),1)\right) = O(|C^{\eta}_{\TT}|),$$
where $|C^{\eta}_{\TT}|$ is the cumulative volume of al cuboids which constitute $C^{\eta}_{\TT}$. Next, observe that
$|C^{\eta}_{\TT}| \to |\TT|$ as $\eta\to 0$. It implies that the contribution by the boundary terms to the variance is $o(1)$,
and thus can be ignored, which eventually means that the asymptotic Gaussianity can be proved by following the same methods as sketched out in this paper,  when considering a $d$-dimensional compact, convex, 
symmetric\footnote{By a symmetric set, we mean that every chord passing through the origin must be bisected
at the origin.} subset of $\real^d$, as parameter space $\TT$.\\

\noindent {\scshape Joint convergence of the various LKCs:}

We note here that using similar ideas, one can prove the multivariate case for different values of the threshold $u$. One of the important questions to look forward to, is the joint distribution of various LKCs evaluated at a fixed threshold.
Although we believe the joint convergence can be proven, getting meaningful estimates
on limiting covariances is likely to be challenging.

Nevertheless, it is worth noticing that the computations done to obtain Theorem \ref{clt-LLk-R} allow to get, in a 
straightforward way, a multivariate CLT for the EPCs $\displaystyle \left(\LL_0(A_u(f;\TT_{U^*_i})), i=1,\cdots,n \right)$:

\begin{cor}[CLT for multivariate $\displaystyle \left(\LL_0(A_u(f;\TT_{U^*_i})), i=1,\cdots,n \right)$] \label{cor:clt-LL0-TUi} ~

\noindent Under Hypothesis $(H1)$ to $(H3)$,  we have, for any $U^*_i\in\Gr(d,k),\,i=1,..,n$, 
\begin{equation*}
 \!\!\left( \LL^{\#}_0(u, \TT;U^*_i) , i=1,\cdots,n\right)^t\; \to\; N\left(0,\Sigma_{0,n}(u)\right),\,\,\,\text{as } \TT\to\real^d,
\end{equation*}
with $\displaystyle \Sigma_{0,n}(u)=\left( \sigma_{0,U^*_i,U^*_j} (u) \right)_{1\le i,j \le n}$ the limiting covariance matrix, $\sigma_{0,U^*_i,U^*_j}(u)$ being the limit, as $\TT\to\real^d$, of the following
\begin{eqnarray*}
&&\frac1{|\TT_{U^*_i}|^{1/2}|\TT_{U^*_j}|^{1/2}}
\sum_{q=1}^{\infty}\sum_{|\un|=q}\sum_{|\um|=q} \!\! c(\un,u,k,\Lambda) c(\um,u,k,\Lambda) \\
&&\qquad \times
\sum_{\scriptsize\begin{array}{c} d_{ij}\ge 0 \\ \sum_{i}d_{ij} = n_j \\ \sum_{j}d_{ij} = m_i\end{array}}
\!\!\!\!\!\!\!\!\!\! \un! \, \um ! \underset{x\in\partial_k \TT_{U^*_i\,}}{\int}\underset{w\in\partial_k (x-\TT_{U^*_j})}{\int} \prod_{1\le i,j\le D} \frac{\gamma^{d_{ij}}_{ij}(w)}{d_{ij}!} \,dw \,dx.
\end{eqnarray*}
Note that the variances are finite and positive. 
\end{cor}

Notice that those finite dimensional distributions given in Corollary \ref{cor:clt-LL0-TUi} might help to obtain the CLT for general LKCs in an alternative way. Indeed, if we may ensure the tightness, then applying the Hadwiger formula \eqref{hadwiger} allows to conclude the CLT of $\LL_k(A_u(f;\TT))$. 
Nevertheless, proving the tightness on such a space is still an open problem.



\section{Acknowledgments} 
Both authors kindly acknowledge the financial support received from IFCAM 
(Indo-French Center for Applied Mathematics) to work on this project in India 
(TIFR - CAM, Bangalore) and in France (ESSEC Business school, Paris) in 
2014 and 2015. This result has been presented at EVA conference (invited `RARE' session) in June 2015.\\
This study has also received the support from the European 
Union's Seventh Framework Programme for research, technological development 
and demonstration under grant agreement no 318984 - RARE, and from the 
{\it Airbus Foundation Chair on Mathematics of Complex Systems} at TIFR-CAM, Bangalore.\\
Note that another study on the same topic (\cite{mue}), has been worked out in parallel providing the same result.

\end{document}